\documentstyle[12pt]{article}
\newtheorem{thm}{Theorem}[section]
\newtheorem{lem}[thm]{Lemma}

\newcommand{\be}{\begin{equation}}
\newcommand{\ee}{\end{equation}}
\newcommand{\bes}{\begin{equation*}}
\newcommand{\ees}{\end{equation*}}

\newcommand{\bea}{\begin{eqnarray}}
\newcommand{\eea}{\end{eqnarray}}
\newcommand{\ba}{\begin{array}}
\newcommand{\ea}{\end{array}}
\newcommand{\bc}{\begin{center}}
\newcommand{\ec}{\end{center}}

\def\s{\sigma}
\def\e{{\bf 1}\!\!{\rm I}}
\def\l{\lambda}
\def\z{\eta}
\def\i{\varepsilon}
\def\t{\tau}

\def\N{\mathbf{N}}
\def\R{\mathbf{R}}
\def\a{\alpha}

\def\m{\mu}
\def\n{\nu}

\def\g{\gamma}

\begin{document}
\begin{center}
{\Large {\bf On stability properties of positive contractions of
$L^1$-spaces accosiated with finite von Neumann algebras}}
\footnote{The work supported by NATO-TUBITAK PC-B programm}\\[4mm]
{\bf Farrukh Mukhamedov}\footnote{Current address: Dipartamento de
Fisica, Universidade de Aveiro, Campus Universitário de Santiago,
3810-193 Aveiro,
Portugal, e-mail: {\tt far75m@yandex.ru}, {\tt farruh@fis.ua.pt} }\\[1mm]

{\it Department of Mechanics and Mathematics\\
National University of Uzbekistan\\
Vuzgorodok, 700174, Tashkent, Uzbekistan}
\\[3mm]
{\bf Hasan Akin}\footnote{e-mail: {\tt hasanakin69@hotmail.com}},
 \ \  {\bf Seyit
Temir}\footnote{e-mail: {\tt seyittemir38@yahoo.com}}\\[1mm]

{\it Department of Mathematics,\\
 Arts and Science Faculty, \\
 Harran University, 63200 \\
 \c{S}anliurfa, Turkey}\\
\end{center}

\begin{abstract}
In the paper we extent the notion of Dobrushin coefficient of
ergodicity  for positive contractions defined on $L^1$-space
associated with finite von Neumann algebra, and in terms of this
coefficient we prove stability results for $L^1$-contractions.
\end{abstract}


\section{Introduction}

It is known (see \cite{K}) that the investigations of asymptotical
behaviors of Markov operators on commutative $L^1$-spaces are very
important. On the other hand, these investigations are related
with several notions of mixing (i.e. weak mixing, mixing,
completely mixing etc) of $L^1$-contractions of a measure space.
Relations between these notions present great interest (see for
example, \cite{BLRT},\cite{BKLM}).  But in  those investigations
it is essentially used that $L^1$-spaces possess the lattice
property. Therefore it is natural to consider Markov operators on
some partially ordered Banach spaces, which are not lattices. One
class of such spaces consists of $L^1$-spaces associated with von
Neumann algebras. It should be noted that this class of Banach
spaces possesses a property of strongly normal cones (see
\cite{EW1}). In \cite{EW1},\cite{EW2},\cite{S} certain
asymptotical properties of Markov semigroups on non-commutative
$L^1$-spaces were studied.

In this paper we are going to study  uniformly (resp. strongly)
asymptotically stable contractions of $L^1$-spaces associated with
finite von Neumann algebras in terms of the Dobrushin
coefficients. The paper is organized as follows. Section 2
contains some preliminary facts and definitions. In Section 3 we
introduce Dobrushin coefficient of ergodicity of
$L^1$-contraction. Using this notion we prove uniformly
asymptotical stability criterion for stochastic operators, which
is a non-commutative analog the Bartoszek's result (see \cite{B}).
Further in section 4 we give an analog of Akcoglu-Sucheston
theorem (see \cite{AS}) for non-commutative $L^1$-spaces. We hope
that this result enables to study subsequential ergodic theorems
in a non-commutative setting (see \cite{CL},\cite{LM}). In the
final section 5 using the results of the previous section we prove
strongly asymptotical stability criterion for positive
$L^1$-contractions. We note that our results are not valid when
von Neumann algebra is semi-finite.

\section{Preliminaries}

Throughout the paper  $M$ would be a von Neumann algebra with the
unit $\e$ and let $\t$ be a faithful  normal finite trace on $M$.
Recall that an element $x\in M$ is called {\it self-adjoint} if
$x=x^*$. The set of all self-adjoint elements is denoted by
$M_{sa}$. By $M_*$ we denote a pre-dual space to $M$.  An element
$p \in M_{sa}$ is called a projector if $p^{2}=p$. Let $\nabla$ be
the set of projectors: $\nabla$ forms a logic. For $p \in \nabla$
we set $p^{\perp}=\e-p$ (see for more definitions
\cite{BR},\cite{T}).

The map $\|\cdot\|_{1}:M\rightarrow [0, \ \infty)$ defined by the
formula $\|x\|_{1}=\tau(|x|)$ is a norm (see \cite{N}). The
completion of $M$ with respect to the norm $\|\cdot\|_{1}$ is
denoted by $L^{1}(M,\t)$. It is known \cite{N} that the spaces
$L^{1}(M,\tau)$  and $M_*$ are isometrically isomorphic, therefore
they can be identified. Further we will use this fact without
noting.

\begin{thm}\label{Lp}\cite{N} The space $L^{1}(M,\tau)$ coincides with the set
$$
L^{1}=\{ x=\int^{\infty}_{-\infty}\lambda
de_{\lambda}:\int^{\infty}_{-\infty}|\lambda| d \tau
(e_{\lambda})< \infty \}.
$$
Moreover,
$$
\|x\|_{1}=\int^{\infty}_{-\infty}|\lambda |d \tau (e_{\lambda}).
$$
Besides,  if $x,y\in L^1(M,\t)$ such that $x\geq 0, y\geq 0$ and
$x\cdot y=0$ then $\|x+y\|_1=\|x\|_1+\|y\|_1$.
\end{thm}

It is known \cite{N} that the equality \be\label{L1}
L^1(M,\t)=L^1(M_{sa},\t)+iL^1(M_{sa},\t) \ee is valid. Note that
$L^1(M_{sa},\t)$ is a pre-dual to $M_{sa}$.

Let $T:L^1(M,\t)\to L^1(M,\t)$ be linear bounded operator.  We say
that a linear operator $T$ is  {\it positive} is $Tx\geq 0$
whenever $x\geq 0$.  A linear operator $T$ is said to be a {\it
contraction} if $\|T(x)\|_1\leq \|x\|_1$ for all $x\in
L^1(M_{sa},\t)$. A positive operator $T$ is called {\it
stochastic} if $\t(Tx)=\t(x)$, $x\geq 0$.  It is clear that any
stochastic operator is a contraction. For given $y\in
L^1(M_{sa},\t)$ and $z\in M_{sa}$ define a linear  operator
$T_{y,z}:L^1(M_{sa},\t)\to L^1(M_{sa},\t)$ as follows
$$
T_{y,z}(x)=\t(xz)y
$$
and extend it to $L^1(M,\t)$ as $T_{y,z}x=T_{y,z}x_1+iT_{y,z}x_2$,
where $x=x_1+ix_2$, $x_1,x_2\in L^1(M_{sa},\t)$.

Put $T_y:=T_{y,\e}$. A linear operator $T:L^1(M,\t)\to L^1(M,\t)$
is called {\it uniformly (resp. strongly) asymptotically stable}
if there exist elements $y\in L^1(M_{sa},\t)$ and $z\in M_{sa}$
such that
$$
\lim_{n\to\infty}\|T^n-T_{y,z}\|=0
$$
(resp. for every $x\in L^1(M,\t)$
$$
\lim_{n\to\infty}\|T^nx-T_{y,z}x\|_1=0.)
$$

\section{Uniformly asymptotically stable contractions}

Let $M$ be a von Neumann algebra with faithful normal finite trace
$\tau$. Let $L^1(M,\tau)$ be a $L^1$-space associated with $M$.

Let $T:L^1(M,\tau)\to L^1(M,\tau)$ be a linear bounded operator.
Define
$$
X=\{x\in L^1(M_{sa},\tau): \ \t(x)=0\},
$$
\be \label{db} \bar\a(T)=\sup_{x\in X,x\neq
0}\frac{\|Tx\|_1}{\|x\|_1}, \ \ \ \ \a(T)=\|T\|-\bar\a(T). \ee

The magnitude $\a(T)$ is called the Dobrushin coefficient of
ergodicity of $T$.

{\bf Remark 3.1.} We note that in commutative case, the notion of
the Dobrushin coefficient of ergodicity was introduced in
\cite{C},\cite{D},\cite{ZZ}.

We have the following theorem which extends the results of
\cite{C},\cite{ZZ}.

\begin{thm}\label{3.1}  Let $T:L^1(M,\tau)\to L^1(M,\tau)$ be a
linear bounded operator. Then the following inequality holds \be
\label{dob} \|Tx\|_1\leq\bar\a(T)\|x\|_1+\a(T)|\t(x)| \ee for
every $x\in L^1(M_{sa},\t)$. \end{thm}

{\bf Proof} Let assume that $x$ is positive. Then $\|x\|_1=\t(x)$
and we have
$$
\bar\a(T)\|x\|_1+\a(T)|\t(x)|=\bar\a(T)\t(x)+(\|T\|-\bar\a(T))\t(x)=\|T\|\|x\|_1\geq\|Tx\|_1.
$$
 So (\ref{dob}) is valid. If $x\leq 0$ the same argument is used
 to prove (\ref{dob}).  Now let $x\in X$ then (\ref{dob}) easily
 follows from (\ref{db}).

 Suppose that $x$ is not in one of the above three cases. Then
 $x=x^+-x^-$, $\|x^+\|_1\neq 0$, $\|x^-\|_1\neq 0$,$\|x^+\|_1\neq \|x^-\|_1$ (see
 \cite{T}). Let  $\|x^+\|_1>\|x^-\|_1$. Put
 $$
 y=\frac{\|x^-\|_1}{\|x^+\|_1}x^+-x^-, \ \ \ \
 z=\frac{\|x^+\|_1-\|x^-\|_1}{\|x^+\|_1}x^+.
 $$
 Then $x=y+z$ and $\|x\|_1=\|y\|_1+\|z\|_1$, here it has been used Theorem 2.1. It is clear that $y\in X$
 and $z\geq 0 $, therefore the inequality (\ref{dob}) is valid for
 $y$ and $z$. Hence, we get
 $$
 \|Tx\|_1\leq\|Ty\|_1+\|Tz\|_1\leq\bar\a(T)\|y\|_1+\bar\a(T)
\|z\|_1+\a(T)\t(z)=\bar\a(T)\|x\|_1+\a(T)|\t(x)|.
 $$

Before formulating the main result of this section we need some
lemmas.

\begin{lem}\label{3.2}  For every $x,y\in L^1(M_{sa},\t)$ such that $x-y\in X$ there
exist $u,v\in L^1(M_{sa},\t)$, $u,v\geq 0$, $\|u\|_1=\|v\|_1=1$,
such that
$$
x-y=\frac{\|x-y\|_1}{2}(u-v).
$$
\end{lem}

{\bf Proof.} We have $x-y=(x-y)^+-(x-y)^-$. Define
$$
u=\frac{(x-y)^+}{\|(x-y)^+\|_1}, \ \ \
v=\frac{(x-y)^-}{\|(x-y)^-\|_1}.
$$
It is clear that $u,v\geq 0$ and $\|u\|_1=\|v\|_1=1$. Since
$x-y\in X$ implies that
$$
\t(x-y)=\t((x-y)^+)-\t((x-y)^-)=\|(x-y)^+\|_1-\|(x-y)^-\|_1=0
$$
therefore $\|(x-y)^+\|_1=\|(x-y)^-\|_1$. Using this and the fact
$\|x-y\|_1=\|(x-y)^+\|_1+\|(x-y)^-\|_1$ we get
$\|(x-y)^+\|_1=\|x-y\|_1/2$. Consequently, we obtain
$$
u-v=\frac{(x-y)^+}{\|x-y\|_1/2}-\frac{(x-y)^-}{\|x-y\|_1/2}=\frac{2}{\|x-y\|_1}(x-y).
$$

\begin{lem}\label{3.3}  Let $T:L^1(M,\t)\to L^1(M,\t)$ be a stochastic operator. Then \be
\label{db2} \bar\a(T)=\sup\left\{\frac{\|Tu-Tv\|_1}{2}: \ u,v\in
L^1(M_{sa},\t), u,v\geq 0, \|u\|_1=\|v\|_1=1 \right\}. \ee
\end{lem}

{\bf Proof.} For $x\in X$, $x\neq 0$ using Lemma \ref{3.2} we have
\bea
\frac{\|Tx\|_1}{\|x\|_1}&=&\frac{\|T(x^+-x^-)\|_1}{\|x^+-x^-\|_1}\nonumber\\
&=&\frac{\frac{\|x^+-x^-\|_1}{2}\|T(u-v)\|_1}{\|x^+-x^-\|_1}\nonumber\\
&=&\frac{\|Tu-Tv\|_1}{2}.\nonumber \eea The equality (\ref{db})
with the last equality imply the desired one (\ref{db2}).\\

Now we are ready to prove the main result of this section, which
is a non-commutative version Bartoszek's result \cite{B}.

\begin{thm}\label{3.4}  Let $T:L^1(M,\t)\to L^1(M,\t)$ be a
stochastic operator.  The following conditions are equivalent:
\begin{enumerate}
   \item[(i)] there exists $\rho>0$ and $n_0\in\N$ such that
   $\a(T^{n_0})\geq \rho$;
   \item[(ii)]  there exits an element $y\in L^1(M_{sa},\t)$, $y\geq 0$ such
   that
$$
\lim_{n\to\infty}\|T^n-T_{y}\|=0.
$$
\end{enumerate}
\end{thm}

{\bf Proof.} (i) $\Rightarrow$ (ii). Let $\rho>0$ and $n_0\in\N$
such that  $\a(T^{n_0})\geq \rho$, this implies that
$\bar\a(T^{n_0})\leq 1-\rho$. Put $\g=1-\rho$.  For an arbitrary
$\i>0$ choose $k\in\N$ such that $\g^k<\i/2$ and set $K=n_0k$.
Since $T$ is a stochastic operator then  we have $\t(T^nx-T^mx)=0$
for every $x\in L^1(M_{sa},\t)$, $x\geq 0$ and
$n,m\in\N\cup\{0\}$. Hence using (\ref{dob}) we infer \bea
\|T^nx-T^mx\|_1&=&\|T^{n_0}(T^{n-n_0}x-T^{m-n_0}x)\|_1\nonumber \\
&\leq & \g\|T^{n-n_0}x-T^{m-n_0}x\|_1\nonumber \\
&\leq&\g^2\|T^{n-2n_0}x-T^{m-2n_0}x\|_1\nonumber \\
&\leq&\cdots\leq\g^k\|T^{n-K}x-T^{m-K}x\|_1\nonumber \\
&\leq&\g^k(\|T^{n-K}x\|_1+\|T^{m-K}x\|_1)\leq 2\g^k\|x\|_1<\i
 \nonumber
 \eea
for every $x\in L^1(M_{sa},\t),x\geq 0$, $\|x\|_1\leq 1$ and
$n,m\geq K$.

Now in general, keeping in mind (\ref{L1}) for every $x\in
L^1(M,\t)$, $\|x\|_1\leq 1$ we have $x=\sum\limits_{k=1}^4i^kx_k$,
$x_k\geq 0$,$\|x_k\|_1\leq 1$, therefore the last relation implies
that
$$
\|T^nx-T^mx\|_1\leq 4\i.
$$
Consequently, we obtain that $(T^n)_{n\in\N}$ is a Cauchy sequence
with respect to uniform norm. Therefore for $x\in L^1(M,\t),x\geq
0$, $\|x\|_1=1$ the sequence $(T^nx)_{n\in\N}$ converges in the
norm of $L^1(M,\t)$ to some $y\in L^1(M,\t)$. Since
$\|Tx\|_1=\|x\|_1=1$ and $T$ is positive, it follows that $y\geq
0$, $\|y\|_1=1$ and $Ty=y$. Using this we obtain
$$
\|T^nz-y\|_1=\|T^nz-T^ny\|_1\leq\|T^{n-1}z-T^{n-1}y\|_1=\|T^{n-1}z-y\|_1
$$
for every $z\in L^1(M_{sa},\t),z\geq 0$, $\|z\|_1\leq 1$.  Hence
the sequence $(\|T^nz-y\|_1)_{n\in\N}$ is monotonic decreasing. By
means of the inequality
$$
\|T^{mn_0}z-y\|_1\leq 2\g^m  \ \ \ \textrm{for every} \ \ m\in\N
$$
we infer that the sequence $(T^nz)_{n\in\N}$ converges to $y$ in
the norm topology of $L^1(M_{sa},\t)$.

If $z\in L^1(M_{sa},\t),z\geq 0$, $\|z\|_1\neq 0$ then taking into
account that
$$
T^nz=\|z\|_1T\bigg(\frac{z}{\|z\|_1}\bigg)=\t(z)T\bigg(\frac{z}{\|z\|_1}\bigg)
$$
we conclude that $T^nz\to\t(z)y$ as $n\to\infty$, since
$T\bigg(\frac{z}{\|z\|_1}\bigg)$ norm converges to $y$.

If $z\in L^1(M_{sa},\t)$, then $z=z^+-z^-$, therefore
$$
T^nz^+\to\t(z^+)y \ \ \textrm{and} \ \ \  T^nz^-\to\t(z^-)y \ \ \
\textrm{as} \ \ n\to\infty.
$$
So $T^nz$ converges to $T_yz$ for every $z\in L^1(M_{sa},\t)$.

In general, if $z\in L^1(M,\t)$, then $z=z_1+iz_2$, where
$z_1,z_2\in L^1(M_{sa},\t)$, hence
$$
T^nz=T^nz_1+iT^nz_2\to\t(z_1)y+i\t(z_2)y=\t(z)y \\ \ \ \textrm{as}
\ \ n\to\infty.
$$
Thus $T^nz$ converges to $T_yz$ for every $z\in L^1(M,\t)$.
 Since $(T^n)_{n\in\N}$ is a Cauchy sequence in
the uniform operator topology it follows that
$$
\lim_{n\to\infty}\|T^n-T_y\|=0.
$$

(ii)$\Rightarrow$ (i). Let $y$ be the element of $L^1(M_{sa},\t)$
defined at (ii). Let $\z\in(0,1/4)$ be given a fixed number. Then
(ii) implies that there is a number $n_0\in\N$ such that
$\|T^n-T_y\|<\z$ for every $n\geq n_0$. Since $Ty=y$ we get that
\be\label{in}
\|T^{n_0}u-T^{n_0}v\|_1\leq\|T^{n_0}u-y\|_1+\|T^{n_0}v-y\|_1<2\z,
\ee for every $u,v\in L^1(M_{sa},\t)$, $u,v\geq 0$,
$\|u\|_1=\|v\|_1=1$.

Hence, using Lemma \ref{3.3} (see (\ref{db2})) we obtain
$\bar\a(T^{n_0})\leq 2\z$ which yields that $\a(T^{n_0})\geq
1-2\z$. The proof is complete.

\section{Completely mixing and smoothing contractions}

In this section we introduce completely mixing  and smoothing
conditions for $L^1$-contractions of non-commutative
$L^1(M,\t)$-space. These notions will be used in next section.

Let $T:L^1(M,\tau)\to L^1(M,\tau)$ be a linear contraction. Define
\be \label{mix}
\bar\rho(T)=\sup\left\{\lim_{n\to\infty}\frac{\|T^n(u-v)\|_1}{\|u-v\|_1}
: \ \ u,v\in L^1(M_{sa},\t), u,v\geq 0, \|u\|_1=\|v\|_1\right\}
\ee and $\rho(T)=\lim\limits_{n\to\infty}\|T^n\|-\bar\rho(T)$.

The magnitude $\rho(T)$ is called the asymptotic  Dobrushin
coefficient of ergodicity of $T$. If $\bar\rho(T)=0$ then $T$ is
called {\it completely mixing}. Note that certain properties  of
completely mixing quantum dynamical systems have been studied in
\cite{AP}.

Using the same argument as in the proof of Theorem \ref{3.4} one
can prove the following

\begin{thm}\label{4.1}  Let $T:L^1(M,\tau)\to L^1(M,\tau)$ be a
linear contraction. Then the following inequality holds \be
\label{mix1}
\lim_{n\to\infty}\|T^nx\|_1\leq\bar\rho(T)\|x\|_1+\rho(T)|\t(x)|
\ee for every $x\in L^1(M_{sa},\t)$. \end{thm}

Using this Theorem we can prove the following

\begin{thm}\label{4.2}  If $T$ is a stochastic operator then
$\bar\rho(T)=0$ or 1. \end{thm}

{\bf Proof} From (\ref{mix}) one can easily see that
$0\leq\bar\rho(T)\leq 1$. Now suppose that $\bar\rho(T)<1$. This
means that there is a number $\g\geq 0$ such that
$\bar\rho(T)\leq\g<1$. Let $x\in X$, $x\neq 0$. It follows that
$$
\lim_{n\to\infty}\|T^nx\|_1\leq\bar\rho(T)\|x\|_1\leq\g\|x\|_1,
$$
therefore, there is a number $n_1\in\N$ such that
$\|T^{n_1}x\|_1\leq\g\|x\|_1$. If $T^{n_1}x=0$ then
$\lim\limits_{n\to\infty}\|T^nx\|_1=0$. If $T^{n_1}x\neq 0$ then
$\t(T^{n_1}x)=\t(x)=0$ since $T$ is stochastic. Thus by means of
(\ref{mix1}) we get
$$
\lim_{n\to\infty}\|T^{n+n_1}x\|_1\leq\bar\rho(T)\|T^{n_1}x\|_1\leq\g\|T^{n_1}x\|_1\leq\g^2\|x\|_1.
$$
It follows that there exists $n_2>n_1$ such that
$\|T^{n_2}x\|_1\leq\g^2\|x\|_1$. Continuing in this way, and if
$T^nx\neq 0$ for every $n\in\N$ then we can find a strictly
increasing sequence $(n_k)$ such that
$\|T^{n_k}x\|_1\leq\g^k\|x\|_1$ for every $k\in\N$. Since $T$ is a
contraction we conclude $\|T^nx\|_1\to 0$ as $n\to\infty$, which
implies that $\bar\rho(T)=0$. \\

Let $T$ be a positive contraction of $L^1(M,\t)$, and let $x\in
L^1(M,\t)$ be such that $x\geq 0$, $x\neq 0$. We say that $T$ is
{\it smoothing} with respect to(w.r.t.) $x$ if for every $\i>0$
there exist $\delta>0$ and $n_0\in\N$ such that $\t(pT^nx)<\i$ for
every $p\in\nabla$ such that $\t(p)<\delta$ and for every $n\geq
n_0$. A commutative counterpart of this notion was introduced in
\cite{ZZ},\cite{KT}. The following result has been proved in
\cite{MTA}, for the sake of completeness we will prove it.

\begin{thm}\label{as-n} Let $T:L^1(M,\t)\to L^1(M,\t)$ be a
positive contraction. Assume that  there is a positive element
$y\in L^1(M,\t)$ such that $T$ is smoothing w.r.t. $y$. Then
$\lim\limits_{n\to\infty}\|T^ny\|_1=0$ or there is a non zero
positive $z\in L^1(M,\t)$ such that $Tz=z$.
\end{thm}

{\bf Proof.} The contractivity of $T$ implies that the limit
$$
\lim_{\n\to\infty}\|T^ny\|_1=\a
$$
exists. Assume that $\a\neq 0$. Define $\l:M_{sa}\to \R$ by
$$
\l(x)=L((\t(xT^ny)_{n\in\N}))
$$
for every $x\in M_{sa}$, here $L$ means a Banach limit (see,
\cite{K}). We have
$$
\l(\e)=L((\t(T^nx)_{n\in\N}))=\lim_{\n\to\infty}\|T^nx\|_1=\a\neq
0,
$$
therefore $\l\neq 0$. Besides, $\l$ is a positive functional,
since for positive element $x\in M_{sa}$,$x\geq 0$ we have
$$
\t(xT^ny)=\t(x^{1/2}T^nyx^{1/2})\geq 0,
$$
for every $n\in\N$.

For arbitrary $x\in M$, we have $x=x_1+ix_2$ and define $\l$ by
$$
\l(x)=\l(x_1)+i\l(x_2).
$$

Let $T^{**}$ be the second dual of $T$, i.e. $T^{**}:M^{**}\to
M^{**}$. The functional $\l$ is $T^{**}$-invariant. Indeed,
\bea\label{inv} (T^{**}\l)(x)&=&<x,T^{**}\l> \nonumber \\
&=&<T^*x,\l>=\nonumber \\
&=&L((\t(T^nyT^*x)_{n\in\N}))\nonumber \\
&=&L((\t(xT^{n+1}y)_{n\in\N}))\nonumber
\\
&=&L((\t(xT^ny)_{n\in\N}))=\l(z).\nonumber \eea

Let $\l=\l_n+\l_s$ be the Takesaki's decomposition (see \cite{T})
of $\l$ onto normal and singular components. Since  $T$ is normal
and $T^{**}\l=\l$, so using the idea of \cite{J} it can be proved
the equality $T^{**}\l_n=\l_n$. Now we will show that $\l_n$ is
nonzero. Consider a measure $\m:=\l{\mid}_\nabla$. It is clear
that $\m$ is an additive measure on $\nabla$. Let us prove that it
is $\s$-additive. To this end, it is enough to show that
$\m(p_k)\to 0$ whenever $p_{k+1}\leq p_k$ and $p_k\searrow 0$,
$p_k\in\nabla$.

Let $\i>0$. From $p_n\searrow 0$ we infer that $\t(p_n)\to 0$ as
$n\to\infty$. It follows that there exists $k_{\i}\in\N$ such that
$\t(p_k)<\i$ for all $k\geq k_\i$. Since $T$ is smoothing w.r.t.
$y$ we obtain
$$
\t(p_kT^ny)<\i, \ \ \ \ \forall k\geq k_\i,
$$
for every $n\geq n_0$. From a property of Banach limit we get
$$
\l(p_k)=L((\t(p_kT^ny)_{n\in\N})<\i \ \ \ \ \textrm{for every} \ \
k\geq k_\i,
$$
which implies $\m(p_k)\to 0$ as $k\to\infty$. This means that the
restriction of $\l_n$ on $\nabla$ coincides with $\m$. Since
$$
\t(p^{\perp}T^ny)>\t(T^ny)-\i\geq\inf\|T^ny\|_1-\i=\a-\i
$$
as $\i$ has been arbitrary, so $\a-\i>0$, and hence
$\m(p^{\perp})>0$ for all $p\in\nabla$ such that $\t(p)<\delta$.
Therefore $\m\neq 0$, and consequently, $\l_n\neq 0$.

From this we infer that there exists a positive element $z\in
L^1(M,\t)$ such that
$$
\l_n(x)=\t(zx), \ \ \ \forall x\in M.
$$
The last equality and $T^{**}\l_n=\l_n$ yield that  \bea
\t(zx)=<x,T^{**}\l_n>=<T^*x,\l_n>=\t(zT^*x)=\t(Tzx) \nonumber \eea
for every $x\in M$, which implies that $Tz=z$.\\

{\bf Remark 4.1} The proved Theorem \ref{as-n} is a
non-commutative analog of Akcoglu and Sucheston result \cite{AS}.
But they used weak convergence instead of smoothing. In fact, the
smoothing is less restrictive than the one they used, since if a
sequence $T^nx$, $x\geq 0$ weakly converges then it is a weak
pre-compact set, hence according to Theorem III.5.4 \cite{T} we
infer that $T$ is smoothing with respect to $x$.

Using the proved Theorem in \cite{MTA} we have proved a
non-commutative analog of \cite{KS} which indicates a relation
between mixing and completely mixing conditions.

{\bf Remark 4.2.} It should be noted that  Theorem \ref{as-n} is
not valid if a von Neumann algebra is semi-finite. Indeed, let
$B(\ell_2)$ be the algebra of all linear bounded operators on
Hilbert space $\ell_2$. Let $\{\phi_n\}$,$n\in\N$ be a standard
basis of $\ell_2$, i.e.
$$
\phi_n=(\underbrace{0,\cdots,0,1}_n,0\cdots).
$$
The matrix units of $B(\ell_2)$ can be defined by
$$
e_{ij}(\xi)=(\xi,\phi_i)\phi_j, \ \ \ \xi\in \ell_2, \ i,j\in\N.
$$
A trace on $B(\ell_2)$ is defined by
$$
\t(x)=\sum_{k=1}^\infty(x\phi_k,\phi_k).
$$
By $\ell_{\infty}$ we denote a maximal commutative subalgebra
generated by elements $\{e_{ii}:\  i\in\N\}$. Let $E:
B(\ell_2)\to\ell_{\infty}$ be the canonical conditional
expectation (see \cite{T}). Define a map
$s:\ell_{\infty}\to\ell_{\infty}$ as follows: for every element
$a\in\ell_{\infty}$, $a=\sum\limits_{k=1}^\infty a_ke_{kk}$ put
$$
s(a)=\sum_{k=1}^\infty a_{k}e_{k+1,k+1}.
$$

Define $T:B(\ell_2)\to B(\ell_2)$ as $T(x)=s(E(x))$, $x\in
B(\ell_2)$. It is clear that $T$ is positive and
$\t(T(x))\leq\t(x)$ for every $x\in L^1(B(\ell_2),\t)\cap
B(\ell_2)$, $x\geq 0$. Hence, $T$ is a positive $L^1$-contraction.
But for this $T$ there is no non zero $x$ such that $Tx=x$.
Moreover, for every $y\in L^1(B(\ell_2),\t)$ we have
$\lim\limits_{n\to\infty}\|T^ny\|_1\neq 0$.

\section{Strongly asymptotical stable contractions}

In this section we give a criterion  for strong asymptotically
stability of contractions in terms of  complete mixing.

Now we are ready to prove a criterion on strong asymptotical
stability.

\begin{thm}\label{5.1} Let $T:L^1(M,\t)\to L^1(M,\t)$ be a
positive contraction. The following conditions are equivalent:
\begin{enumerate}
   \item[(i)]  $T$ is completely mixing and smoothing w.r.t. some
   $h\in L^1(M,\t)$, $h\geq 0$;
   \item[(ii)]  there exits an element $y\in L^1(M,\t)$, $y\geq 0$ such
   that for every $x\in L^1(M,\t)$
$$
\lim_{n\to\infty}\|T^nx-T_{y}x\|_1=0.
$$
\end{enumerate}
\end{thm}

{\bf Proof.} (i)$\Rightarrow$(ii). Let $h\in L^1(M,\t)$, $h\geq
0$, $h\neq 0$ be such that $T$ is smoothing w.r.t. $h$. Without
loss of generality we may assume that $\|h\|_1=1$. By Theorem
\ref{as-n} there are only two possibilities:

(a) \  $\lim\limits_{n\to\infty}\|T^nh\|_1=0$;

(b) \ there exists $y\in L^1(M,\t)$, $y\geq 0$, $y\neq 0$ such
that $Ty=y$.

If we are in the situation (a), then for every  $x\in
L^1(M,\t)$,$x\geq 0$,$\|x\|_1=1$ using complete mixing one gets
$$
\lim_{n\to\infty}\|T^nx\|_1\leq
\lim_{n\to\infty}\|T^nx-T^nh\|_1+\lim_{n\to\infty}\|T^nh\|_1=0.
$$

Let $x\in L^1(M,\t)$, then we have $x=\sum\limits_{k=1}^4i^kx_k$,
where $x_k\geq 0$. Using this from the last relation one gets that
$T^n$ converges strongly to $T_0$.

Let us consider the second situation (b). In this case we may
assume that $\|y\|_1=1$. Completely mixing condition for $T$
implies that
$$
\lim_{n\to\infty}\|T^nx-y\|_1=0
$$
for every $x\in L^1(M,\t)$, $x\geq 0$, $\|x\|_1=1$. The similar
arguments used towards the end of the proof of Theorem \ref{3.4}
show the desired relation holds.

(ii)$\Rightarrow$ (i). Let $g\in X$ then $T^ng$ norm converges to
$\t(g)y=0$, and hence  $T$ is completely mixing.

Let $x\in L^1(M,\t)$, $x\geq 0$, $\|x\|_1=1$, then the sequence
$(T^nx)$ norm converges to $y$. So according to Remark 4.1 we find
that $T$ is smoothing w.r.t. $x$.\\

{\bf Acknowledgements.}  The first named author (F.M.) thanks
TUBITAK for providing financial support and Harran University  for
kind hospitality and providing all facilities. The work is also
partially supported by Grants $\Phi$-1.1.2, $\Phi$.2.1.56 of CST
of the Republic of
Uzbekistan.\\


\end{document}